\documentclass[a4paper]{amsart}

\usepackage{amsmath,amsthm,amssymb,amscd}
\usepackage[arrow,matrix]{xy}
\usepackage[dvips]{graphicx}
\usepackage{enumerate}

\theoremstyle{plain}
\numberwithin{equation}{section}
\newtheorem{thm}{Theorem}[section]
\newtheorem{prop}[thm]{Proposition}

\newtheorem{lem}[thm]{Lemma}
\theoremstyle{definition}
\newtheorem{dfn}[thm]{Definition}
\newtheorem{exm}[thm]{Example}

\def\rank{\mathop{\mathrm{rank}}\nolimits}
\def\dim{\mathop{\mathrm{dim}}\nolimits}
\def\Lie{\mathop{\mathrm{Lie}}\nolimits}

\def\kutorsor{E_{\sigma}\to\Gamma(\sigma)^{\mathrm{gp}} \backslash D_{\sigma}}
\def\gsds{\Gamma(\sigma)^{\mathrm{gp}} \backslash D_{\sigma}}
\def\gs{\Gamma(\sigma)^{\mathrm{gp}}}
\def\gm{\mathbb{G}_{m}}
\def\es{E_{\sigma}}

\def\c{\mathbb{C}}
\def\q{\mathbb{Q}}
\def\r{\mathbb{R}}
\def\z{\mathbb{Z}}

\def\im{\mathop{\mathrm{Im}}\nolimits}
\def\bs{\backslash}
\def\hom{\mathop{\mathrm{Hom}}\nolimits}
\def\Aut{\mathop{\mathrm{Aut}}\nolimits}
\def\torus{\mathrm{torus}}
\def\toric{\mathrm{toric}}
\def\spec{\mathop{\mathrm{Spec}}\nolimits}

\newcommand{\mf}[1]{{\mathfrak{#1}}}
\newcommand{\mb}[1]{{\mathbf{#1}}}
\newcommand{\bb}[1]{{\mathbb{#1}}}
\newcommand{\mca}[1]{{\mathcal{#1}}}
\title[moduli of log Hodge structures]{On the boundary of the moduli spaces \\of log Hodge structures: \\triviality of the torsor}
\author[T.~Hayama]{Tatsuki HAYAMA}
\thanks{Supported by Grant-in-Aid for JSPS Fellows from Japan Society for the Promotion of Science.}
\address{Graduate School of Science, Osaka University, Osaka 560-0043, Japan}
\email{t-hayama@cr.math.sci.osaka-u.ac.jp}

\begin{document}
\maketitle
\begin{abstract}
In this paper we will study the moduli spaces of log Hodge structures introduced by Kato-Usui. This moduli space is a partial compactification of a discrete quotient of a period domain. We treat  the following 2 cases: (A) the case where the period domain is Hermitian symmetric, (B) the case where the Hodge structures are of the mirror quintic type. Especially we study a property of the torsor.
\end{abstract}
\section{Introduction}
Let $\varphi :(\Delta^*)^n\to \Gamma\backslash D$ be a period map arising from a variation of Hodge structures with unipotent monodromies on the $n$-fold product of the punctured disk. Here $\Gamma$ is the image of $\pi_1((\Delta^*)^n)\cong\z^n$ by the monodromy representation, i.e., $\Gamma$ is a free $\z$-module. By Schmid's nilpotent orbit theorem \cite{s}, the behavior of the period map around the origin is approximated by a ``nilpotent orbit". Then we add the set of nilpotent orbits to $\Gamma\backslash D$ as boundary points and extend the period map satisfying the following diagram:
\begin{align*}
\xymatrix{
(\Delta )^n\ar@{}[d]|{\rotatebox{270}{$ \supset $}}\ar@{->}[r]&\Gamma\backslash D\cup \{\text{nilpotent orbits}\}\ar@{}[d]|{\rotatebox{270}{$ \supset $}}
\\
(\Delta^*)^n\ar@{->}[r]^{\varphi }&\Gamma\backslash D.\\
}
\end{align*}
Here $\Gamma\backslash D$ is an analytic space and $\varphi $ is an analytic morphism. But, except in some cases, we have no way to endow the upper right one with an analytic structure. 

In \cite{ku} Kato-Usui endow it with a geometric structure as a ``logarithmic manifold" and they treat the above diagram as a diagram in the category of logarithmic manifolds. Moreover they define ``polarized logarithmic Hodge structures" (PLH for abbr.) and they show that the upper right one is the moduli space of PLH. Our result is for the geometric structure of the moduli spaces of PLH in the  following 2 cases:
\begin{enumerate}
\item[(A)] $D$ is a Hermitian symmetric space.
\item[(B)] Polarized Hodge structures (PH for abbr.) are of weight $w=3$ and Hodge number $h^{p,q}=1\;(p+q=3, p,q\geq 0)$ and logarithms of the monodromy transformations are of type $N_{\alpha}$ in \cite[\S 12.3]{ku} (or type ${\rm II}_1$ in \cite{GGK}).
\end{enumerate}
The case (A) corresponds to degenerations of algebraic curves or K3 surfaces. This case is classical and well-known. 
The case (B) corresponds to degenerations of certain Calabi-Yau three-folds, for instance those occurring in the mirror quintic family. This case is studied recently.
For example, Griffiths et al (\cite{GGK}) describe ``N\'eron models'' for VHS of this type. Usui (\cite{u}) shows a logarithmic Torelli theorem for this quintic mirror family. 

\subsection*{Construction of a moduli space of PLH}
To explain our result, we describe Kato-Usui's construction of the moduli space of PLH roughly (\S 3 for detail). Steps of the construction are given as follows:
\begin{enumerate}
\item[{\bf Step 1.}] Define the nilpotent cone $\sigma$ and the set $D_{\sigma}$ of nilpotent orbits.
\item[{\bf Step 2.}] Define the toric variety $\toric_{\sigma}$ and the space $E_{\sigma}$. 
\item[{\bf Step 3.}] Define the map $E_{\sigma} \to\Gamma\backslash D_{\sigma}$ and endow $\Gamma\backslash D_{\sigma}$ with a geometric structure by this map. 
\end{enumerate}

Firstly, we fix a point $s_0\in(\Delta^*)^n$ and let $(H_{s_0},F_{s_0},\langle\; ,\;\rangle_{s_0})$ be the corresponding polarized Hodge structure (PH for abbr.).
The period domain $D$ is a homogeneous space for the real Lie group $G=\Aut{(H_{s_0,\r},\langle\; ,\;\rangle_{s_0})}$ and also an open $G$-orbit in the flag manifold $\check{D}$ (\S 2 for detail). 

{\bf Step 1}: Take the logarithms $N_1,\ldots ,N_n$ of the monodromy transformations and make the cone $\sigma$ in $\mf{g}$ generated by them, which is called a {\it nilpotent cone}. 
By Schmid's nilpotent orbit theorem, there exists the limiting Hodge filtration $F_{\infty}$. We call the orbit $\exp{(\sigma_{\c})}F_{\infty}$ in $\check{D}$ a {\it nilpotent orbit}. 
$D_{\sigma}$ is the set of all nilpotent orbits generated by $\sigma$. 

{\bf Step 2}: Take the monoid $\Gamma(\sigma):=\Gamma\cap\exp{(\sigma)}$.
It determines the toric variety $\toric_{\sigma}:=\spec(\c[\Gamma(\sigma)^{\vee}])_{\mathrm{an}}$. We define the subspace $E_{\sigma}$ of $\toric_{\sigma}\times\check{D}$ in \S 3.4.

{\bf Step 3}: We define the map $E_{\sigma} \to\Gamma\backslash D_{\sigma}$ in \S 3.5 (here $\Gamma=\Gamma(\sigma)^{\text{gp}}$). By {\cite[Theorem A]{ku}}, $E_{\sigma}$ and $\Gamma\backslash D_{\sigma}$ are logarithmic manifolds. Moreover the map is a $\sigma_{\c}$-torsor in the category of logarithmic manifolds, i.e., there exists a proper and free $\sigma_{\c}$-action on $E_{\sigma}$ and $E_{\sigma} \to\Gamma\backslash D_{\sigma}$ is isomorphic to  $E_{\sigma} \to\sigma_{\c}\backslash E_{\sigma}$ in the category of logarithmic manifolds.

\subsection*{Main result} 
Our main result is for properties of the torsor $E_{\sigma} \to\Gamma\backslash D_{\sigma}$.

\begin{itemize}\label{main1}
\item In the case (A), the torsor is trivial. (Theorem \ref{mainthm})
\item In the case (B), the torsor is non-trivial. (Proposition \ref{mainthm2})
\end{itemize}

In the case (A), $\Gamma\bs D_{\sigma}$ is just a toroidal partial compactification of $\Gamma\bs D$ introduced by \cite{amrt}.
To show the triviality, we review the theory of bounded symmetric domains in \S 4.
Realization of $D$ as the Siegel domain of the 3rd kind (\ref{S3isom}) is a key of the proof.
This induces the triviality of $B(\sigma)\to\bb{B}(\sigma)$ (Lemma \ref{triv}).
By the triviality of this torsor, we show the triviality of $E_{\sigma} \to\Gamma\backslash D_{\sigma}$.
We also describe a simple example (Example \ref{exm}).
 
In the case (B), $D$ is not a Hermitian symmetric domain, i.e., isotropy subgroups are not maximally compact.
We fix a point $F_0\in D$ and take a maximally compact subgroup $K$ as in (\ref{maxcompact}).
Then $K F_0$ is a compact subvariety of $D$.
Existence of such a variety (of positive dimension) is a distinction between the cases where the period domain is Hermitian symmetric and otherwise.
The compact subvariety plays an important role in the proof.

\subsection*{Acknowledgment}

A part of main results is from the author's master thesis.
The author is grateful to Professors Sampei Usui, Christian Schnell for their valuable advice and warm encouragement. 
The author is also grateful to the referee for his careful reading and valuable suggestions and
comments on presentations.

\section{Polarized Hodge structures and period domains} 
We recall the definition of polarized Hodge structures and of period domains.
A {\it Hodge structure} of weight $w$ and of Hodge type $(h^{p,q})$ is a pair $(H,F)$ consisting of a free $\z$-module of rank $\sum_{p,q}h^{p,q}$ and of a decreasing filtration on $H_{\c}:=H\otimes \c$ satisfying the following conditions:
\begin{enumerate}
\item[(H1)] $\dim_{\c} F^p=\sum_{r\geq p}h^{r,w-r}\quad \text{for all $p$.}$\label{check1}
\item[(H2)] $H_{\c}=\bigoplus _{p+q} H^{p,q}\quad(H^{p,q}:=F^p\cap \overline{F^{w-p}}).$
\end{enumerate}

A {\it polarization} $\langle\; ,\;\rangle$ for a Hodge structure $(H,F)$ of weight $w$ is a non-degenerate bilinear form on $H_{\q}:=H\otimes \q$, symmetric if $w$ is even and skew-symmetric if $w$ is odd, satisfying the following conditions:
\begin{enumerate}
\item[(P1)] $\langle F^p , F^q \rangle=0\quad \text{for $p+q>w$.}$\label{check2}
\item[(P2)] The Hermitian form 
$$H_{\c}\times H_{\c}\to\c,\quad (x,y)\mapsto \langle C_{F}(x),\bar{y}\rangle$$
is positive definite.
\end{enumerate}
Here $\langle\; ,\;\rangle$ is regarded as the natural extension to $\c$-bilinear form and $C_{F}$ is the Weil operator, which is defined by $C_{F}(x):=(\sqrt{-1})^{p-q}x$ for $x\in H^{p,q}$.
 
We fix a polarized Hodge structure $(H_{0},F_{0},\langle\; ,\;\rangle_0)$ of weight $w$ and of Hodge type $(h^{p,q})$.
We define the set of all Hodge structures of this type
$$D:=\left\{\begin{array}{l|l}F&\begin{array}{r}(H_{0} , F ,\langle\; ,\;\rangle_0)\text{ is a polarized Hodge structure}\\ \text{ of weight $w$ and of Hodge type $(h^{p,q})$}\end{array}\end{array}\right\}.$$
$D$ is called a {\it period domain}.
Moreover, we have the flag manifold 
$$\check{D}:=\left\{\begin{array}{l|l}F&\begin{array}{r}(H_0 , F ,\langle\; ,\;\rangle_0)\text{ satisfy the conditions}\\ \text{ (H1), (H2) and (P1)}\end{array}\end{array}\right\}.$$
$\check{D}$ is called the {\it compact dual of $D$}.
$D$ is contained in $\check{D}$ as an open subset.
$D$ and $\check{D}$ are homogeneous spaces under the natural actions of $G$ and $G_{\c}$ respectively, where $G:=\Aut{(H_{0,\r},\langle\; ,\;\rangle_0)}$ and $G_{\c}$ is the complexification of $G$.
$G$ is a classical group such that
\begin{align*}
G\cong\begin{cases}
Sp(h,\r)/\{\pm 1\}& \text{if $w$ is odd,}\\
SO(h_{\text{odd}},h_{\text{even}})& \text{if $w$ is even,}\\
\end{cases}
\end{align*}
where $2h=\rank{H_0}$, $h_{\text{odd}}=\sum_{\text{p:odd}}h^{p,q}$ and $h_{\text{even}}=\sum_{\text{p:even}}h^{p,q}$.
The isotropy subgroup of $G$ at $F_0$ is isomorphic to
\begin{align*}
\begin{cases}
\prod_{p\leq m}U(h^{p,q})& \text{if $w=2m+1$,}\\
\prod_{p<m}U(h^{p,q})\times SO(h^{m,m})& \text{if $w=2m$.}\\
\end{cases}
\end{align*}
They are compact subgroups of $G$ but not maximal compact in general.
$D$ is a Hermitian symmetric domain if and only if the isotropy subgroup is a maximally compact subgroup, i.e., one of the following is satisfied:
\begin{enumerate}
\item $w=2m+1, \;h^{p,q}=0 \text{ unless }p=m+1,m.$
\item $w=2m, \;h^{p,q}=1 \text{ for }p=m+1,m-1,\;h^{m,m}\text{ is arbitary, }h^{p,q}=0\text{ otherwise. }$
\item $w=2m, \;h^{p,q}=1 \text{ for }p=m+a,m+a-1,m-a,m-a+1\text{ for some }a\geq 2,\;h^{p,q}=0\text{ otherwise. }$
\end{enumerate} 
In the case (1), $D$ is a Hermitian symmetric domain of type {\rm III}.
In the case (2) or (3), an irreducible component of $D$ is a Hermitian symmetric domain of type {\rm IV}.
\section{Moduli spaces of polarized log Hodge structures}
In this section, we review some basic facts in \cite{ku}.
Firstly, we introduce $D_{\Sigma}$, a set of nilpotent orbits associated with a fan $\Sigma$ consisting of nilpotent cones.
Secondly, for some subgroup $\Gamma$ in $G_{\z}$, we endow $\Gamma\bs D_{\Sigma}$ with a geometric structure.
Finally we see some fundamental properties of $\Gamma\bs D_{\Sigma}$, which are among the main results of \cite{ku}.  
\subsection{Nilpotent orbits}
A {\it nilpotent cone} $\sigma$ is a strongly convex and finitely generated rational polyhedral cone in  $\mf{g}:=\Lie{G}$ whose generators are nilpotent and commute with each other. 
For $A=\r,\c$, we denote by $\sigma_{A}$ the $A$-linear span of $\sigma$ in $\mf{g}_A$.
\begin{dfn}\label{nilp}
Let $\sigma=\sum_{j=1}^n\r_{\geq 0}N_j$ be a nilpotent cone and $F\in\check{D}$. 
$\exp{(\sigma_{\c})}F\subset\check{D}$ is a {\it $\sigma$-nilpotent orbit} if it satisfies following conditions:
\begin{enumerate}
\item $\exp{(\sum_j iy_jN_j)}F\in D$ for all $y_j\gg 0$.
\item $NF^p\subset F^{p-1}$ for all $p\in\z$ and for all $N\in \sigma$.
\end{enumerate}
\end{dfn}
The condition (2) says the map $\c^n\to \check{D}$ given by $(z_j)\mapsto\sum_j\exp(z_jN_j)F$ is horizontal.
Let $\Sigma$ be a fan consisting of nilpotent cones. We define the set of nilpotent orbits
$$D_{\Sigma}:=\{(\sigma,Z)|\; \sigma\in\Sigma,\; Z\;\text{is a}\;\sigma \text{-nilpotent orbit}\}.$$
For a nilpotent cone $\sigma$, the set of faces of $\sigma$ is a fan, and we abbreviate $D_{\{\text{faces of }\sigma\}}$ as $D_{\sigma}$.
\subsection{Subgroups in $G_{\z}$ which is compatible with a fan}
Let $\Gamma$ be a subgroup of $G_{\z}$ and $\Sigma$ a fan of nilpotent cones.
We say $\Gamma$ is {\it compatible} with $\Sigma$ if 
$$\mathrm{Ad}(\gamma)(\sigma)\in\Sigma$$
for all $\gamma\in\Gamma$ and for all $\sigma\in\Sigma$.
Then $\Gamma$ acts on $D_{\Sigma}$ if $\Gamma$ is compatible with $\Sigma$.
Moreover we say $\Gamma$ is {\it strongly compatible} with $\Sigma$ if it is compatible with $\Sigma$ and for all $\sigma\in \Sigma$ there exists $\gamma_1,\ldots,\gamma_n\in \Gamma(\sigma):=\Gamma\cap\exp{(\sigma)}$ such that
$$\sigma=\sum_{j}\r_{\geq 0}\log{(\gamma_j)}.$$

\subsection{Varieties $\toric_{\sigma}$ and $\torus_{\sigma}$}
Let $\Sigma$ be a fan and $\Gamma$ a subgroup of $G_{\z}$ which is strongly compatible with $\Sigma$.
We have toric varieties associated with the monoid $\Gamma(\sigma)$ such that
\begin{align*}
&\toric_{\sigma}:=\spec(\c[\Gamma(\sigma)^{\vee}])_{\mathrm{an}}\cong \hom{(\Gamma(\sigma)^{\vee},\c)},\\
&\torus_{\sigma}:=\spec(\c[{\Gamma(\sigma)}^{\vee \mathrm{gp}} ])_\mathrm{an}\cong\hom{({\Gamma(\sigma)}^{\vee \mathrm{gp}},\gm)}\cong \gm\otimes {\Gamma(\sigma)}^{\mathrm{gp}},
\end{align*}
where $\c$ is regarded as a semigroup via multiplication and above homomorphisms are of semigroups. 
As in \cite[\S 2.1]{F}, we choose the distinguished point 
$$x_{\tau}:\Gamma(\sigma)^{\vee}\to\c;\quad u\mapsto \begin{cases}1&\text{if }u\in \Gamma(\tau)^{\perp },\\0&\text{otherwise,}\end{cases}$$
for a face $\tau$ of $\sigma$.
Then $\toric_{\sigma}$ can be decomposed as 
$$\toric_{\sigma}=\bigsqcup_{\tau\prec \sigma}(\torus_{\sigma}\cdot x_{\tau}).$$
For $q\in\toric_{\sigma}$, there exists $\sigma(q)\prec\sigma$ such that $q\in\torus_{\sigma}\cdot x_{\sigma(q)}$.
By a surjective homomorphism
\begin{equation}\label{emap}
\mb{e}:\sigma_{\c}\rightarrow \torus_{\sigma}\cong\gm\otimes {\Gamma(\sigma)}^{\mathrm{gp}};\;w\log{(\gamma)}\mapsto \exp{(2\pi\sqrt{-1}w)}\otimes \gamma ,
\end{equation}
$q$ can be written as
\begin{equation}\label{orbit decomp}
q=\mb{e}(z)\cdot x_{\sigma (q)}.
\end{equation}
Here $\ker{(\mb{e})}=\log{(\Gamma(\sigma)^{\text{gp}})}$ and $z$ is determined uniquely modulo $\log{(\gs)}+\sigma(q)_{\c}$.

\subsection{Spaces $\check{E}_{\sigma}$ and $E_{\sigma}$}
We define the analytic space $\check{E}_{\sigma}:=\toric_{\sigma}\times \check{D}$ and endow $\check{E}_{\sigma}$ with the logarithmic structure $M_{\es}$ by the inverse image of canonical logarithmic structure on $\toric_{\sigma}$ (cf.\ \cite{k}).
Define the subspace of $\check{E}_{\sigma}$
\begin{equation*}
E_{\sigma}:=\{(q,F)\in\check{E}_{\sigma}\;|\;\exp{(\sigma(q)_{\c})}\exp{(z)}F\text{ is $\sigma (q)$-nilpotent orbit}\}
\end{equation*}
where $z$ is an element such that $q=\mb{e}(z)\cdot x_{\sigma (q)}$.
The set $E_{\sigma}$ is well-defined.
The topology of $\es$ is the ``strong topology" in $\check{E}_{\sigma}$, which is defined in \cite[\S 3.1]{ku}, and $\mca{O}_{\es}$(resp.\ $M_{E_{\sigma}}$) is the inverse image of $\mca{O}_{\check{E}_{\sigma}}$(resp.\ $M_{\check{E}_{\sigma}}$).
Then $E_{\sigma}$ is a logarithmic local ringed space.
Note that $E_{\sigma}$ is {\it not} an analytic space in general.

\subsection{The structure of $\Gamma\backslash D_{\Sigma}$}
We define the canonical map
\begin{align*}\label{es}
&E_{\sigma}\to \Gamma(\sigma)^{\mathrm{gp}}\backslash D_{\sigma},\\
&(q,F)\mapsto (\sigma(q),\exp{(\sigma(q)_{\c})}\exp{(z)}F)\mod{\Gamma(\sigma)^{\mathrm{gp}}},
\end{align*}
where $q=\mb{e}(z)\cdot x_{\sigma (q)}$ as in (\ref{orbit decomp}).
This map is well-defined.
We endow $\Gamma\backslash D_{\Sigma}$ with the strongest topology for which the composite maps $\pi_{\sigma}:\es\to\gsds\to\Gamma\backslash D_{\Sigma}$ are continuous for all $\sigma\in\Sigma$.
We endow $\Gamma\backslash D_{\Sigma}$ with $\mca{O}_{\Gamma\backslash D_{\Sigma}}$ (resp.\ $M_{\Gamma\backslash D_{\Sigma}}$) as follows:
\begin{align*}&\mathcal{O}_{\Gamma\backslash D_{\Sigma}}(U)\,(\text{resp.\ }M_{\Gamma\backslash D_{\Sigma}}(U))\\
&:=\{\mathrm{map}\;f:U\rightarrow \c|\;f\circ  \pi_{\sigma}\in\mathcal{O}_{E_{\sigma}}(\pi_{\sigma}^{-1}(U))\,(\text{resp.\ }M_{E_{\sigma}}(\pi_{\sigma}^{-1}(U)))\;({}^{\forall}\sigma\in\Sigma)\}\end{align*}
for any open set $U$ of $\Gamma\backslash D_{\Sigma}$.
As for $E_{\sigma}$, note that $\Gamma\backslash D_{\Sigma}$ is a logarithmic local ringed space but is {\it not} an analytic space in general.
Kato-Usui introduce ``logarithmic manifolds" as generalized analytic spaces (cf.\ \cite[\S 3.5]{ku}) and they show the following geometric properties of $\Gamma\backslash  D_{\Sigma}$:
\begin{thm}[{\cite[Theorem A]{ku}}]
Let $\Sigma$ be a fan of nilpotent cones and let $\Gamma$ be a subgroup of $G_{\z}$ which is strongly compatible with $\Sigma$. 
Then we have
\begin{enumerate}
\item $E_{\sigma}$ is a logarithmic manifold.
\item If $\Gamma$ is neat $($i.e., the subgroup of $\gm$ generated by all the eigenvalues of all $\gamma\in\Gamma$ is torsion free$)$, $\Gamma\backslash D_{\Sigma}$ is also a logarithmic manifold.
\item Let $\sigma\in \Sigma$ and define the action of $\sigma_{\c}$ on $E_{\sigma}$ over $\gsds$ by
$$a\cdot(q,F):=(\mb{e}(a)q,\exp{(-a)}F)\quad(a\in\sigma_{\c},\;(q,F)\in \es).$$
Then $\kutorsor$ is a $\sigma_{\c}$-torsor in the category of logarithmic manifold.
\item $\gsds\to\Gamma\bs D_{\Sigma}$ is open and locally an isomorphism of logarithmic manifold.
\end{enumerate}
\end{thm}
In \cite[\S 2.4]{ku} , Kato-Usui introduce ``polarized log Hodge structures" and they show that $\Gamma\backslash D_{\Sigma}$ is a fine moduli space of polarized log Hodge structures if $\Gamma$ is neat(\cite[Theorem B]{ku}).

\section{The structure of bounded symmetric domains}
In this section we recall some basic facts on Hermitian symmetric domains (for more detail, see \cite[III]{amrt} , \cite[appendix]{n}).
We define Satake boundary components, and show that a Hermitian symmetric domain is a family of tube domains parametrized by a vector bundle over a Satake boundary component.
This domain is called a Siegel domain of third kind. 

\subsection{Satake boundary components}
Let $D$ be a Hermitian symmetric domain. 
Then $\Aut{(D)}$ is a real Lie group and the identity component $G$ of $\Aut{(D)}$ acts on $D$ transitively.
We fix a base point $o\in D$.
The isotropy subgroup $K$ at $o$ is a maximally compact subgroup of $G$.
Let $s_o$ be a symmetry at $o$ and let
\begin{align*}
&\mf{g}:=\Lie{(G)},\quad \mf{k}:=\Lie{(K)},\\
&\mf{p}:=\text{the subspace of }\mf{g} \text{ where } s_o=-\mathrm{Id}. 
\end{align*}
Then we have a Cartan decomposition
$$\mf{g}=\mf{k}\oplus\mf{p}.$$
$\mf{p}$ is isomorphic to the tangent space to $D$ at $o$.
Let $J$ be a complex structure on $\mf{p}$ and let
\begin{align*}
&\mf{p}_+:=\text{the $\sqrt{-1}$-eigenspace for $J$ in $\mf{p}_{\c}$},\\
&\mf{p}_-:=\text{the $-\sqrt{-1}$-eigenspace for $J$ in $\mf{p}_{\c}$}. 
\end{align*}
Here $\mf{p}_+$ and $\mf{p_-}$ are abelian subalgebras of $\mf{g}_{\c}$.
Then we have the Harish-Chandra embedding map $D\to \mf{p}_+$ whose image is a bounded domain.
\begin{dfn}
A {\it Satake boundary component} of $D$ is an equivalence class in $\overline{D}$, the topological closure of $D$ in $\mf{p}_{+}$, under the equivalence relation generated by $x\sim y$ if there exists a holomorphic map
$$\lambda :\{z\in\c|\;|z|<1\}\longrightarrow \mf{p}_+$$
such that $\im(\lambda)\subset\overline{D}$ and $x,y\in\im(\lambda)$.
\end{dfn}
It is known that Satake boundary components are also Hermitian symmetric domains.
Let $S$ be a Satake boundary component.
We define
\begin{align*}
N(S)&:=\{ g\in G|\; gS=S\},\\
W(S)&:=\;\text{the unipotent radical of $N(S)$},\\
U(S)&:=\; \text{the center of $W(S)$}.
\end{align*}
These groups have the following properties:
\begin{prop}\label{groups}
\begin{enumerate}
\item $N(S)$ acts on $D$ transitively.
\item There exists an abelian Lie subalgebra $\underline{v}(S)\subset \mf{g}$ such that $\Lie{W(S)}=\underline{v}(S)+\Lie{U(S)}$.
\item $W(S)/U(S)$ is an abelian Lie group which is isomorphic to $V(S):=\exp{\underline{v}(S)}$.
\end{enumerate}
\end{prop}
$S$ is called {\it rational} if $N(S)$ is defined over $\q$.
If $S$ is rational then $V(S)$ and $U(S)$ are also defined over $\q$. 
\subsection{Siegel domain of third kind}
We define a subspace of $\check{D}$
$$D(S):=U(S)_{\c}\cdot D=\bigcup_{g\in U(S)_{\c}}g\cdot D$$
where $U(S)_{\c}:=U(S)\otimes \c$.
By Proposition \ref{groups}(1) and the fact that $U(S)$ is a normal subgroup, $N(S)U(S)_{\c}$ acts on $D(S)$ transitively.
We choose the base point $o_{S}$ as in \cite[\S 4.2]{amrt} .
The isotropy subgroup $I$ of $G$ at $o_{S}$ is contained in $N(S)$.
Then we have a map
$$\Psi _{S}:D(S)\cong N(S)U(S)_{\c}/I\to N(S)U(S)_{\c}/N(S)\cong U(S)$$
where the last isomorphism takes imaginary part.
By \cite[\S 4.2 Theorem 1]{amrt} , we have an open homogeneous self adjoint cone $C(S)\subset U(S)$ such that
$\Psi _{S}^{-1}(C(S))=D$.
\begin{thm}\label{S3}
\begin{enumerate}
\item $U(S)_{\c}$ acts freely on $D(S)$. $D(S)\to U(S)_{\c}\bs D(S)$ is trivial principal homogeneous bundle.
\item $V(S)$ acts freely on $U(S)_{\c}\bs D(S)$. $V(S)\bs (U(S)_{\c}\bs D(S))\cong S$ and the quotient map $D(S)\to S$ is a {\it complex} vector bundle (although $V(S)$ is {\it real}). Moreover it is trivial.
\item By (1) and (2), we have a trivialization 
\begin{equation}\label{S3ds}
D(S)\cong S\times \c^k \times U(S)_{\c}.
\end{equation}
In this product representation,we have
$$\Psi_S(x,y,z)=\im{z}-h_x(y,y)$$
where $h_x$ is a real-bilinear quadratic form $\c^k\times\c^k\to U(S)$ depending real-analytically on $x$. 
\end{enumerate}
\end{thm}
Thus we have 
\begin{equation}\label{S3isom}
D\cong\left\{\begin{array}{l|l} (x,y,z)\in S\times \c^{(g-k)k} \times U(S)_{\c}&\im z\in C(S)+h_x(y,y)\end{array}\right\}.
\end{equation}
\section{Main result}
\subsection{The case: $\kutorsor$ is trivial}
In this subsection, we assume that $D$ is a period domain and also a Hermitian symmetric domain. 
The purpose is to show Theorem \ref{mainthm}.
Main theorem for the case where $D$ is upper half plane is described in Example \ref{exm}.

Let $S$ be a Satake rational boundary component of $D$ and $\sigma$ a nilpotent cone included in $\Lie{(U(s))}$. 
Firstly we show the triviality of the torsor for such a cone.
We set 
$$B(\sigma):=\exp{(\sigma_{\c})}\cdot D\subset\check{D},\quad\bb{B}(\sigma):=\exp{(\sigma_{\c})}\backslash B(\sigma).$$
Here $B(\sigma),\bb{B}(\sigma)$ are defined by Carlson-Cattani-Kaplan (\cite{cck}) from the point of view of mixed Hodge theory. 
\begin{lem}\label{triv}
$B(\sigma)\to\bb{B}(\sigma)$ is a trivial principal bundle with fiber $\exp{(\sigma_{\c})}$.
\end{lem}
\begin{proof}
$\exp{(\sigma_{\c})}$ is a sub Lie group of $U(S)_{\c}$ and $U(S)_{\c}$ is abelian.
By Theorem \ref{S3}(1), $D(S)\to\exp{(\sigma_{\c})}\bs D(S)$ is a trivial principal bundle .
Furthermore the following diagram is commutative:
\begin{align*}
\xymatrix{
D(S)\ar@{}[d]|{\rotatebox{270}{$ \supset $}}\ar@{->}[r]&
\exp{(\sigma_{\c})}\bs D(S)\ar@{}[d]|{\rotatebox{270}{$ \supset $}}\\
B(\sigma)\ar@{->}[r]&\bb{B}(\sigma)\\
}
\end{align*}
Then $B(\sigma)\to\bb{B}(\sigma)$ is trivial.
\end{proof}
Now we describe a trivialization of $B(\sigma)$ over $\bb{B}(\sigma)$ explicitly. 
Take a complementary sub Lie group $Z_{\sigma}$ of $\exp{(\sigma_{\c})}$ in $U(S)_{\c}$.
By (\ref{S3ds}), we have a decomposition of $D(S)$ as $D(S)\cong S\times\c^k\times Z_{\sigma}\times \exp{(\sigma_{\c})}$.
Here 
\begin{equation}\label{dstriv}
\exp{(\sigma_{\c})}\bs D(S)\cong S\times\c^k\times Z_{\sigma}
\end{equation}
and the decomposition of $D(S)$ makes a trivialization of $D(S)$ over $\exp{(\sigma_{\c})}\bs D(S)$.
Via (\ref{dstriv}), we have $Y\subset S\times\c^k\times Z_{\sigma}$ satisfying the following commutative diagram:
\begin{align*}
\xymatrix{
\exp{(\sigma_{\c})}\bs D(S)\ar@{}[d]|{\rotatebox{270}{$ \supset $}}\ar@{}[r]|{\rotatebox{0}{$ \cong $}}&
S\times\c^k\times Z_{\sigma}\ar@{}[d]|{\rotatebox{270}{$ \supset $}}\\
\bb{B}(\sigma)\ar@{}[r]|{\rotatebox{0}{$ \cong $}}&Y\\
}
\end{align*}
In fact, $Y$ is the image of $D$ via the projection $D(S)\to S\times\c^k\times Z_{\sigma}$.
Then $B(\sigma)\cong \exp{(\sigma_{\c})}\times Y$ is a trivialization of $B(\sigma)$ over $\bb{B}(\sigma)$.

Let $\Gamma$ be a subgroup of $G_{\z}$ which is strongly compatible with $\sigma$.
Let us think about a quotient trivial bundle $\gs\backslash B(\sigma)\to\bb{B}(\sigma)$.
Its fiber is the quotient of $\exp{(\sigma_{\c})}$ by the lattice $\Gamma(\sigma)^{\text{gp}}$.
Since $\exp{(\sigma_{\c})}$ is a unipotent and abelian Lie group, $\sigma_{\c}\cong\exp{(\sigma_{\c})}$.
Via this isomorphism, the lattice action on $\exp{(\sigma_{\c})}$ is equivalent to the lattice action on $\sigma_{\c}$ by $\log{(\gs)}$. 
Then the fiber is isomorphic to
$$\sigma_{\c}/\log{(\gs)}=\sigma_{\c}/\ker{(\mb{e})}\cong \torus_{\sigma}\quad\text{by (\ref{emap}).}$$
By the canonical torus action, $\gs\backslash B(\sigma)\to\bb{B}(\sigma)$ is also a principal $\torus_{\sigma}$-bundle whose trivialization is given by
\begin{equation*}
\torus_{\sigma}\times Y\stackrel{\sim}{\to}\gs\bs B(\sigma);\quad (\mb{e}(z),F)\mapsto \exp{(z)}F\mod{\gs},
\end{equation*}
where we regard $Y$ as a subset of $\check{D}$ via (\ref{S3ds}).
By the definition of $E_{\sigma}$, we have 
\begin{lem} The following diagram is commutative:
\begin{align*}
\xymatrix{
\gs\backslash B(\sigma)\ar@{}[d]|{\rotatebox{270}{$ \supset $}}\ar@{}[r]|{\rotatebox{0}{$ \cong $}}&
\torus_{\sigma}\times Y\ar@{}[d]|{\rotatebox{270}{$ \supset $}}\\
\gs\bs D\ar@{}[r]|{\rotatebox{0}{$ \cong $}}&(\torus_{\sigma}\times Y)\cap E_{\sigma}\\
}
\end{align*}
\end{lem}
By the torus embedding $\torus_{\sigma}\hookrightarrow \toric_{\sigma}$, we can construct the associated bundle 
$$(\gs\backslash B(\sigma))_{\sigma}:=(\gs\backslash B(\sigma))\times_{\torus_{\sigma}}\toric_{\sigma},$$
and define
\begin{align*}
(\gs\backslash D)_{\sigma}:=\text{ the interior of the closure of }\gs\backslash D\\
\text{ in }(\gs\backslash B(\sigma))_{\sigma}.
\end{align*}
This is a toroidal partial compactification {\it associated with $\sigma$}.
\begin{lem}\label{last}
The following diagram is commutative:
\begin{align*}
\xymatrix{
(\gs\backslash B(\sigma))_{\sigma}\ar@{}[d]|{\rotatebox{270}{$ \supset $}}\ar@{}[r]|{\rotatebox{0}{$ \cong $}}&
\toric_{\sigma}\times Y\ar@{}[d]|{\rotatebox{270}{$ \supset $}}\\
(\gs\bs D)_{\sigma}\ar@{}[r]|{\rotatebox{0}{$ \cong $}}&(\toric_{\sigma}\times Y)\cap E_{\sigma}\\
}
\end{align*}
\end{lem}
\begin{proof}
For $(q,F)\in\toric_{\sigma}\times Y$, $(q,F)\in E_{\sigma}$ if and only if $\exp{(\sigma(q)_{\c})}\exp{(z)}F$ is a $\sigma{(q)}$-nilpotent orbit, where $q=\exp{(z)}x_{\sigma(q)}$ as in (\ref{orbit decomp}).
Since $D$ is Hermitian symmetric, the horizontal tangent bundle of $D$ coincides with the tangent bundle of $D$.
Then the condition of Definition \ref{nilp}(2) is trivially satisfied.
Let $\{N_i\}$ be a set of rational nilpotent elements generating $\sigma(q)$.
$(q,F)\in E_{\sigma}$ if and only if $\exp{\left(\Sigma_j\; y_jN_j\right)}\exp{(z)}F\in D$, i.e.,
\begin{equation*}
\left(\mb{e}\left(\Sigma_j\; y_jN_j+z\right),F\right)\in (\torus_{\sigma}\times Y)\cap E_{\sigma}
\end{equation*}
for all $y_j$ such that $\im{(y_j)}\gg 0$.
In $\toric_{\sigma}$,
$$\lim_{\im{(y_j)}\to\infty}\mb{e}\left(\Sigma_j\; y_jN_j+z\right)=\mb{e}(z)x_{\sigma(q)}=q.$$
Then $(q,F)$ is in the interior of the closure of $(\torus_{\sigma}\times Y)\cap E_{\sigma}$.
\end{proof}
The map $(\toric_{\sigma}\times Y)\cap E_{\sigma}\hookrightarrow E_{\sigma}\to \gsds$ is bijective.
By \cite[8.2.7]{ku}, $E_{\sigma}$ is an open set in $\check{E}_{\sigma}$.
Then $E_{\sigma}\to\gsds$ is a $\sigma_{\c}$-torsor in the category of analytic spaces and
$$(\toric_{\sigma}\times Y)\cap E_{\sigma}\cong \sigma_{\c}\backslash E_{\sigma}\cong \gsds.$$
Thus $(\toric_{\sigma}\times Y)\cap E_{\sigma}\hookrightarrow E_{\sigma}$ gives a section of the torsor $\kutorsor$, i.e. $\kutorsor$ is trivial.

Next we show that $\kutorsor$ is trivial for {\it all} nilpotent cones $\sigma$.
Let $\Gamma=G_{\z}$.
By \cite[{\rm II}]{amrt} , there exists $\Gamma$-admissible collection of fans $\Sigma=\{\,\Sigma(S)\}_S$ where $\Sigma(S)$ is a fan in $\overline{C(S)}$ for every Satake rational boundary component $S$.
Taking logarithm, we identify $\Sigma$ with the collection of fans in $\mf{g}$ which are strongly compatible with $\Gamma$.
We show that $\Sigma$ is large enough to cover all nilpotent cones, i.e., $\Sigma$ is complete fan.

Let $U(S)_{\z}=U(S)\cap\Gamma$.
To obtain $U(S)_{\z}\bs D_{\sigma}$, we should confirm the following proposition:
\begin{prop}
Generators of $Z_{\sigma}$ can be taken in $G_{\z}$.  
\end{prop}
\begin{proof}
$\Gamma(\sigma)^{\text{gp}}$ is saturated in $U(S)_{\z}$, i.e.,
$$\text{If }g\in U(S)_{\z}\text{ and }g^n\in \gs\text{ for some $n\geq 1$, then } g\in \gs.$$ 
Then $U(S)_{\z}/\Gamma(\sigma)^{\text{gp}}$ is a free module and there exists a subgroup $Z_{\sigma,\z}$ in $U(S)_{\z}$ such that $U(S)_{\z}=\Gamma(\sigma)^{\text{gp}}\oplus Z_{\sigma,\z}$.
Hence we have $U(S)_{\c}=\exp{\sigma_{\c}}\oplus (Z_{\sigma,\z}\otimes \c)$.
\end{proof}
Gluing $U(S)_{\z}\backslash D_{\sigma}=Z_{\sigma,\z}\backslash (\gsds)$ for $\sigma\in\Sigma(S)$, we have $U(S)_{\z}\backslash D_{\Sigma(S)}$ as a toroidal partial compactification {\it in the direction $S$}.
Then we obtain a compact variety $\Gamma\backslash D_{\Sigma}$ by \cite{amrt}
(Take quotient by $N(S)_{\z}/ U(S)_{\z}$ and glue neighborhoods of boundaries of $(N(S)_{\z}/ U(S)_{\z})\backslash (U(S)_{\z}\backslash D_{\Sigma(S)})$ with $\Gamma\backslash D$).

On the other hand, we have the following proposition.
\begin{prop}[{\cite[12.6.4]{ku}}]Let $\Gamma$ be a subgroup of $G_{\z}$ and let $\Sigma$ be a fan which is strongly compatible with $\Gamma$. Assume that $\Gamma\backslash D_{\Sigma}$ is compact. Then $\Sigma$ is complete.
\end{prop}
The precise definition of complete fan is given in \cite{ku}.
An important property of a complete fan $\Sigma$ is the following: if there exists $Z\subset\check{D}$ such that $(\sigma ,Z)$ is a nilpotent orbit, then $\sigma\in \Sigma$.
Now $\Gamma$-admissible collection of fans $\Sigma$ is complete since $\Gamma\backslash D_{\Sigma}$ is compact.
It is to say that a nilpotent cone $\sigma$ has a $\sigma$-nilpotent orbit if 
$$\exp{(\sigma)}\subset \overline{C(S)}\subset U(S)$$
 for some Satake boundary component $S$, and $\sigma$ has no $\sigma$-nilpotent orbit ,i.e., $D_{\sigma}=D$, otherwise.
Hence we have
\begin{thm}\label{mainthm}
Let $\sigma$ be a nilpotent cone in $\mf{g}$.
Then $\kutorsor$ is trivial. 
\end{thm}
\begin{exm}\label{exm}
Let $D$ be the upper half plane.
$G=SL(2,\r)$ acts on $D$ by linear fractional transformation.
By the Cayley transformation, $D\cong \Delta:=\{\;z\in\c\,|\,|z|<1\;\}$.
Take the Satake boundary component $S=\{1\}\in \partial  \Delta$.
Then
\begin{align*}&N(S)=\left\{\begin{array}{l|l}\begin{pmatrix}u&v\cr 0&u^{-1}\end{pmatrix}&u\in\r\setminus\{0\},\;v\in \r\end{array}\right\},\\
&W(S)=U(S)=\left\{\begin{array}{l|l}\begin{pmatrix}1&v\cr 0&1\end{pmatrix}&v\in \r\end{array}\right\},\\
&C(S)=\left\{\begin{array}{l|l}\begin{pmatrix}1&v\cr 0&1\end{pmatrix}&v\in \r_{\geq 0}\end{array}\right\}
\end{align*}
(cf.\ \cite{n}).
Take the nilpotent $N=\begin{pmatrix}0&1\cr 0&0\end{pmatrix}$ and the nilpotent cone $\sigma=\r_{\geq 0}N\subset\mf{g}$.
Here $\exp{(\sigma_{\c})}=U(S)_{\c}$.
The compact dual $\check{D}$ and the subspace $B(\sigma)=\exp{(\sigma_{\c})}\cdot D\subset \check{D}$ are described as
$$\check{D}=\c\sqcup \{\infty\}\cong\bb{P}^1,\quad B(\sigma)=\c.$$
$\exp{(\sigma_{\c})}\bs B_{\sigma}=\bb{B}(\sigma)$ is a point and $B(\sigma)\to\bb{B}(\sigma)$ is a trivial principal bundle over $\bb{B}(\sigma)$.
Take $F_0\in B(\sigma)$.
We have a trivialization
$$B(\sigma)=\exp{(\sigma_{\c})}\cdot F_0\cong \exp{(\sigma_{\c})}\times \{F_0\}.$$
For $\Gamma=SL(2,\z)$, 
$$\Gamma(\sigma)^{\text{gp}}=\left\{\begin{array}{l|l}\begin{pmatrix}1&v\cr 0&1\end{pmatrix}&v\in \z\end{array}\right\}=\exp{(\z N)}$$
is a lattice of $\exp{(\sigma_{\c})}$ and $\bb{G}_{m}\cong\exp{(\sigma_{\c})}/\gs$.
Then we have the trivial principal $\bb{G}_{m}$-bundle $\gs\bs B(\sigma)\to\bb{B}(\sigma)$.
By the torus embedding $\bb{G}_{m}\hookrightarrow \c$, we have the trivial associated bundle
\begin{align*}
(\gs\bs B(\sigma))\times_{\gm}\c\cong \c\times\{F_0\}.
\end{align*}

And we have
\begin{align*}
&E_{\sigma}=\left\{\begin{array}{l|l}\left( q,F\right)\in\c\times B(\sigma)&
\begin{array}{l}\exp{((2\pi i)^{-1}\log{(q)}N)}F\in D\text{ if $q\neq 0$}\\
\text{and }F\in B(\sigma )\text{ if $q=0$}
\end{array}
\end{array}\right\}\end{align*}
The map 
\begin{align*}
&(\c\times \{F_0\})\cap E_{\sigma}\hookrightarrow E_{\sigma}\to\gsds ;\\
&(q,F_0)\to\left\{
\begin{array}{ll}
(\{0\},\exp{((2\pi i)^{-1}\log{(q)}N)}F_0) \mod{\gs}& \text{ if $q\neq 0$,}\\
(\sigma,\c) & \text{ if $q= 0$.}\end{array}
\right.
\end{align*}
is an isomorphism by Lemma \ref{last} and we have the following commutative diagram:
\begin{align*}
\xymatrix{
(\c\times \{F_0\})\cap E_{\sigma}\ar@{}[d]|{\rotatebox{90}{$ \cong  $}}\ar@{}[r]|{\rotatebox{0}{\quad\quad $ \subset  $}}&
E_{\sigma}\ar@{->}[dl]\\
\gs\bs D_{\sigma}.\\
}
\end{align*}
Then the torsor $\kutorsor$ has a section, i.e., the torsor is trivial.
\end{exm}
\subsection{The case: $\kutorsor$ is non-trivial}
Let $w=3$, and $h^{p,q}=1\;(p+q=3, p,q\geq 0)$. 
Let $H_0$ be a free module of rank $4$, $\langle\; ,\;\rangle_0$ a non-degenerate alternating bilinear form on $H_0$.
In this case $D\cong Sp(2,\r)/(U(1)\times U(1))$.
Then $D$ is not a Hermitian symmetric space. 
Take $e_1,\ldots ,e_4$ as a symplectic basis for $(H_0,\langle\; ,\;\rangle_0)$, i.e.,
\begin{align*}
\left(\langle e_i,e_j\rangle_0 \right)_{i,j}=
\begin{pmatrix} 0&-I\cr I&0\cr
\end{pmatrix}.
\end{align*}
Define $N \in\mf{g}$ as follows:
\begin{align*}
N(e_3)=e_1,\quad N(e_j)=0\;(j\neq 3).
\end{align*}

\begin{prop}\label{mainthm2}
Let $\sigma=\r_{\geq 0}N$. Then $\kutorsor$ is non-trivial.\end{prop}
\begin{proof}
Define
\begin{align*}
&\left( u_1,\ldots ,u_4 \right):=\frac{1}{\sqrt{2}}(e_1,\ldots ,e_4)
\begin{pmatrix} I&I\cr iI&-iI\cr
\end{pmatrix}\text{, i.e.,}\\
&\left( \langle u_i,u_j\rangle_0 \right)_{i,j}=
\begin{pmatrix} 0&iI\cr -iI&0\cr
\end{pmatrix},\quad (\langle u_i,\overline{u_j}\rangle_0)_{i,j}=
\begin{pmatrix} iI&0\cr 0&-iI\cr
\end{pmatrix}.
\end{align*}
Take $F_w,F_{\infty}\in D\;(w\in\c)$ as follows:
\begin{align*}
&F_w^3=\text{span}_{\c}\{wu_1+u_2\},\quad F_w^2=\text{span}_{\c}\{wu_1+u_2,u_3-wu_4\},\\
&F_{\infty}^3=\text{span}_{\c}\{u_1\},\quad F_{\infty}^2=\text{span}_{\c}\{u_1,u_4\}.
\end{align*}
We have the maximally compact subgroup of $G$ at $F_0$
\begin{align}\label{maxcompact}
K&=\left\{\begin{array}{l|l} g\in G&\langle C_{F_0}(gv),\overline{gw}\rangle_0=\langle C_{F_0}(v),\overline{w}\rangle_0 \quad\text{for }v,w\in H_{\c}\end{array}\right\}\\
&=\left\{\begin{array}{l|l}\begin{pmatrix}
X&0\cr
0&\overline{X}\cr
\end{pmatrix}
&X\in U(2) \end{array}\right\},\nonumber
\end{align}
where matrices in the second equation are expressed with respect to the basis $(u_1,\ldots ,u_4)$.
The $K$-orbit of $F_0$ is given by
\begin{align*}
K\cdot F_{0}=K_{\c}\cdot F_0=\left\{\begin{array}{l|l}F_w&w\in\c\end{array}\right\}\sqcup F_{\infty}\cong \bb{P}^1.
\end{align*}

We assume that $\kutorsor$ is trivial.
Let $\varphi $ be a section of $\kutorsor$.
We define a holomorphic morphism $\Phi :D\to\c$ such that
$$\Phi :D \stackrel{\text{quot.}}{\longrightarrow }\Gamma(\sigma)^{\text{gp}}\backslash D\hookrightarrow \gsds \stackrel{\varphi }{\longrightarrow }E_{\sigma}\stackrel{\text{proj.}}{\longrightarrow }\toric_{\sigma}\cong\c.$$
Since $K F_0\cong \bb{P}^1$, $\Phi |_{K F_{0}}$ is constant.

On the other hand, $(\sigma ,\exp{(\sigma_{\c})}F_0)$ is a nilpotent orbit (it is easy to check the condition of Definition \ref{nilp}).
Then 
\begin{align}\label{limit}
\lim_{x\to\infty}\Phi(\exp{(ixN)}F_0)=0.
\end{align}
Define $N'\in\mf{g}$ as follows:
\begin{align*}
N'(u_3)=u_1,\quad N'(u_j)=0\;(j\neq 3).
\end{align*}
Then we have
\begin{align}
&\exp{(ixN)}F_0=\exp{\left(\frac{x}{2+x}N'\right)}F_0,\label{f0}\\
&F_{\infty}=\exp{\left(\frac{x}{2+x}N'\right)}F_{\infty}\label{finf}
\end{align}
for $x\in\r\setminus\{-2\}$ and 
\begin{align}
\exp(zN')K F_0\subset D \quad \text{ for }|z|<1.\label{ind}
\end{align}
Then $\Phi |_{\exp(zN')K F_{0}}$ is constant for each $|z|<1$, again because $\exp{(zN')}KF_0\cong \mathbb{P}^1$.
Finally we have
\begin{align*}
\begin{array}{rll}
\Phi(\exp{(ixN)}F_0)&=\Phi\left(\exp{\left(\frac{x}{2+x}N'\right)}F_0\right)&\quad\text{ (by (\ref{f0}))}\\
&=\Phi\left(\exp{\left(\frac{x}{2+x}N'\right)}F_{\infty}\right)&\quad\text{ (by (\ref{ind}) and $|\frac{x}{2+x}| <1$)}\\
&=\Phi (F_{\infty})&\quad\text{ (by (\ref{finf}))}
\end{array}
\end{align*}
for $x>-1$.
This contradicts the condition (\ref{limit}), since $\Phi(F_{\infty})\in\torus_{\sigma} $ if $F_{\infty}\in D$.
\end{proof}

\end{document}